\documentclass[12pt,a4paper]{article}

\usepackage[latin2]{inputenc}
\usepackage[T1]{fontenc}
\usepackage[english]{babel}
\usepackage{amsmath}
\usepackage{amssymb}
\usepackage{graphicx}
\usepackage{tikz}
\usepackage{caption, subcaption, float}
\usepackage{verbatim}
\usepackage{dsfont}
\usepackage{enumitem}
\usetikzlibrary{shapes,decorations.pathmorphing}

\newtheorem{thr}{Theorem}
\newtheorem{dfn}{Definition}
\newenvironment{prf}{\noindent\textit{Proof.}\relax}{\vspace{3mm}}%
\newenvironment{rmk}{\noindent\textit{Remark.}\relax}{\vspace{3mm}}%
\newenvironment{step}[1]{\vspace{3mm}\noindent\textbf{Step #1}\relax}{}%
\newtheorem{lmm}{Lemma}
\newtheorem{prp}{Proposition}

\newtheorem{obs}{Observation}
\newtheorem{cnj}{Conjecture}
\newtheorem{prb}{Problem}

\newcommand{\BB}{\mathcal{B}}

\DeclareMathAlphabet\mathbfcal{OMS}{cmsy}{b}{n}

\title{A Characterization of Uniquely Representable Graphs}
\date{\today}
\author{P\'eter G.N. Szab\'o\thanks{This work was supported by the National Research, Development and Innovation Office -- NKFIH,  No. 108947.}\\
Email: szape@cs.bme.hu\\
Department of Computer Science and Information Theory\\
Budapest University of Technology and Economics\\
Budapest, H-1111, Hungary}

\begin{document}

\maketitle

\begin{abstract}
The betweenness structure of a finite metric space $M = (X, d)$ is a pair $\BB(M) = (X,\beta_M)$ where $\beta_M$ is the so-called betweenness relation of $M$ that consists of point triplets $(x, y, z)$ such that $d(x, z) = d(x, y) + d(y, z)$.
The underlying graph of a betweenness structure $\BB = (X,\beta)$ is the simple graph $G(\BB) = (X, E)$ where the edges are pairs of distinct points with no third point between them. A connected graph $G$ is uniquely representable if there exists a unique metric betweenness structure with underlying graph $G$.
It was implied by previous works that trees are uniquely representable. In this paper, we give a characterization of uniquely representable graphs by showing that they are exactly the block graphs. Further, we prove that two related classes of graphs coincide with the class of block graphs and the class of distance-hereditary graphs, respectively. We show that our results hold not only for metric but also for almost-metric betweenness structures.
\end{abstract}


\section{Introduction}\label{S1}

\subsection{Historic Overview}\label{S11}
Metric space is a universal concept with various applications in physics, molecular biology and phylogenetics as well as in different branches of math\-e\-matics including geometry, topology, graph theory, computer science and algebra. Although finite metric spaces are trivial objects from a topological point of view, they have rich combinatorial properties which were investigated from different angles over the last century. The concept of metric betweenness appears to play a central role in the related literature which span from combinatorial geometry to metric graph theory.

An important line of study in metric graph theory was concerned with the axiomatization of different graph classes in terms of geodesic betweenness and more generally, in terms of different path transit functions.
Further, graph classes satisfying certain transit axioms were characterized in different ways, for example, by means of forbidden induced subgraphs. Among the earliest results are Sholander's axiomatic characterization of trees, lattices and partially ordered sets in terms of segments, medians and betweenness \cite{sholander1952trees}. 

The (geodesic) interval function of a connected graph is one of the most fundamental concepts in metric graph theory. It was first extensively studied by Mulder \cite{mulder1980interval}, who introduced the five classical axioms of the interval function and used the interval function to study different classes of graphs.
In a series of papers starting from 1994, Nebesk\'{y} et al. presented a characterization of the interval function of a connected graph in terms of first order transit axioms, each time improving the proof \cite{nebesky1994characterization, nebesky1994characterization-set, nebesky1998characterizing, nebesky2001characterization, nebesky2001interval, mulder2009axiomatic}. His approach was to extend the set of the five classical axioms with additional postulates.

Besides the interval function, two other path transit functions were extensively studied on connected graphs: the induced path function and the all-paths function. In \cite{mulder2008transit}, Mulder introduced the general notion of transit function to unify the three concepts and presented a list of prototype problems that can be uniformly applied to different kinds of transit functions.

In \cite{changat2001all-paths}, Changat et al. characterized the all-paths function of a connected graph by using only first-order transit axioms.
Interestingly, Nebesk\'{y} established in \cite{nebesky2002induced} that such a characterization of the induced path function $J$ of an arbitrary connected graph is impossible. Nevertheless, Changat et al. succeeded in axiomatizing $J$ for certain classes of graphs \cite{changat2003induced, changat2010induced}. The authors also investigated different axioms on $J$ such as the Peano axiom, the Pasch axiom and monotonicity \cite{changat2010induced, changat2004induced, changat2016induced}.
Several well-known graph classes --including trees, block graphs, bipartite graphs, chordal graphs and Ptolemaic graphs-- were characterized by axiomatizing their interval function or induced path function \cite{changat2013forbidden, balakrishnan2015axiomatic, changat2018axiomatic}. For a thorough survey on geodesic and induced path betweenness, see \cite{changat2019betweenness}.

Beside the above mentioned results, there are two other lines of research that are related to the topic of this paper.
Based on the pioneering work of Isbell \cite{isbell1964six} and Buneman \cite{buneman1974note}, Dress et al. developed T-theory, a combinatorial theory of discrete metric spaces, and extensively studied algorithmic and combinatorial aspects of phylogenetic trees \cite{dress1987parsimonious, bandelt1992canonical}. These results have a wide range of applications in evolutionary biology.
Another notable problem that is related to metric betweenness and gained a lot of attention lately is the generalization of the de Bruijn--Erd\H{o}s theorem known from geometry to finite metric spaces, originally conjectured by Chen and Chv\'atal in 2008 \cite{chen2008problems}. The conjecture is still open but has already been proved for several important classes of metric spaces and graph metrics \cite{kantor2013debruijn, chvatal2014debruijn, aboulker2015chen-chvatal, beaudou2015debruijn, aboulker2016lines, aboulker2018new, beaudou2019bisplit, schrader2019debruijn}.

The aim of our research in the most general sense is to better understand the combinatorial properties of the metric betweenness relation. We believe that finite metric spaces are interesting objects on their own right, not only in relation with graphs, posets or Euclidean geometry. Therefore, there are some important differences between our approach and that of metric graph theory.
First, instead of studying metric properties of graphs to draw graph-theoretic conclusions, we focus on the combinatorial aspects of finite metric spaces with the occasional help of graph-theoretic tools and arguments.
Second, we formulate our results in the language of the metric betweenness relation. It is essentially equivalent to the language of transit functions but better suited for describing potentially non-graphic problems.
Third, we do not restrict our attention to betweenness structures obtained in some way from simple graphs but only require them to be metric or almost-metric (see Definition \ref{Dalmmetr}). Results obtained in this framework are general enough to be useful in many different areas of mathematics. Although our characterization does involve graphic betweenness structures, graphicity does not play a distinguished role.

\subsection{Definitions}\label{S12}

A \emph{finite metric space} is a pair $M = (X, d)$ where $X$ is a finite nonempty set and $d$ is a \emph{metric} on $X$, i.e. an $X\times X\rightarrow\mathbb{R}$ function which satisfies the following conditions for all $x, y, z\in X$:
\begin{enumerate}
\item $d(x, y) = 0\Leftrightarrow x = y$ (\emph{identity of indiscernibles})\label{Ems1};
\item $d(x, y) = d(y, x)$ (\emph{symmetry})\label{Ems2};
\item $d(x, z)\leq d(x, y) + d(y, z)$ (\emph{triangle-inequality})\label{Ems3}.
\end{enumerate}
The non-negativity of metric follows from the definition, as for all $x, y\in X$, $0 = d(x, x)\leq d(x, y) + d(y, x) = 2d(x, y)$.
We will refer to the ground set and the metric of a metric space $M$ by $X(M)$ and $d_M$, respectively.
All metric spaces in this paper will be assumed to be \emph{finite} ($|X(M)| <\infty$) if not stated otherwise.

In order to capture the combinatorial properties of metric spaces that are relevant to us, we introduce the following abstraction.
A \emph{betweenness structure} is a pair $\BB = (X,\beta)$ where $X$ is a nonempty finite set and $\beta\subseteq X^3$ is a ternary relation, called the \emph{betweenness relation} of $\BB$.
The fact $(x, y, z)\in\beta$ will be denoted by $(x\ y\ z)_\BB$ or simply by $(x\ y\ z)$ if $\BB$ is clear from the context, and we say that $y$ is \emph{between} $x$ and $z$.

We say that a betweenness structure $\mathcal{C} = (X,\gamma)$ is an \emph{extension} of $\BB$ (in notation, $\BB\preccurlyeq\mathcal{C}$) if $\beta\subseteq\gamma$. The \emph{substructure} of $\BB$ induced by a nonempty subset $Y\subseteq X$ is the betweenness structure $\BB\vert_Y = (Y,\beta\cap Y^3)$.

We believe that stating and proving our results on the abstraction level of betweenness structures makes our arguments clearer and easier to understand, as they are about combinatorial properties of the betweenness relation of metric spaces.

There is a natural way to associate a betweenness structure with a metric space. The \emph{betweenness structure induced by a metric space} $M = (X, d)$ is $\BB(M) = (X,\beta_M)$ where $$\beta_M =\{(x, y, z)\in X^3: d(x, z) = d(x, y) + d(y, z)\}$$ is the \emph{betweenness relation} of $M$. To simplify notations, we will write $(x\ y\ z)_M$ for $(x\ y\ z)_{\BB(M)}$.

The betweenness structure $\BB$ is \emph{metric} if it is induced by some metric space $M = (X, d)$.
A metric betweenness relation satisfies the following elementary properties:\\
\indent for all $x, y, z\in X$,
\begin{enumerate}[label={(P\arabic*)}, ref={(P\arabic*)}]
\item $(x\ x\ z)$\label{Ecoll1};
\item $(x\ y\ z)\Rightarrow (z\ y\ x)$\label{Ecoll2};
\item $(x\ y\ z)\wedge (y\ x\ z)\Rightarrow x = y$\label{Ecoll3};
\end{enumerate}
The \emph{trichotomy} of betweenness follows straight from properties \ref{Ecoll1}--\ref{Ecoll3}: for any three distinct points $x, y, z\in X$, at most one of the relations $(x\ y\ z)$, $(y\ z\ x)$, $(z\ x\ y)$ can hold.

Additionally, for all $x, y, z, w\in X$, we have
\begin{enumerate}[label={(P\arabic*)}, ref={(P\arabic*)}]
\setcounter{enumi}{3}
\item $(x\ y\ z)\wedge (x\ w\ y)\Rightarrow (x\ w\ z)\wedge (w\ y\ z)$.\label{Ecoll4}
\end{enumerate}
This property, that we call the \emph{four relations property},
is the simplest non-trivial property of metric betweennesses.
Notice that properties \ref{Ecoll1}--\ref{Ecoll4} are not sufficient to guarantee that the betweenness structure is metric (think, for example, about the Fano plane, which is proved to be non-metric \cite{chvatal2004sylvester, chen2006sylvester}).

\begin{dfn}\label{Dalmmetr}
A betweenness structure is \emph{almost-metric} if it satisfies properties \ref{Ecoll1}--\ref{Ecoll4}.
\end{dfn}

Betweenness structures of this kind are often called ``pseudometric'' in the related literature, however, we want to avoid confusion with a different meaning of this term, a generalized metric where zero distances are allowed.

Different notions of betweenness were developed by Sholander \cite{sholander1952trees}, Morgana \cite{morgana2002induced} and Burigana \cite{burigana2009tree}. In this context, our concept of almost-metric betweenness can be best viewed as a stronger version of Morgana's betweenness. There, trichotomy and the condition $(x\ y\ z)\wedge (x\ w\ y)\Rightarrow (x\ w\ z)$ is assumed in place of properties \ref{Ecoll3} and \ref{Ecoll4}. To see that our version is indeed a strengthening, consider the set $X=\{w, x, y, z\}$ of four points with betweenness relation $\beta =\{(w, y, z), (z, y, w), (w, x, y), (y, x, w), (w, x, z), (z, x, w)\}\cup\{(u, u, v): u, v\in X\}\cup\{(u, v, v): u, v\in X\}$. The betweenness structure $\BB = (X,\beta)$ is not almost-metric as it violates property \ref{Ecoll4}, but $\beta$ is a betweenness in the sense of Morgana.

Quite interestingly, several properties of finite metric spaces can be generalized to almost-metric betweenness structures without modification. The reader will see that our main results fall in this category as well.

\begin{figure}[t]
\centering
\includegraphics{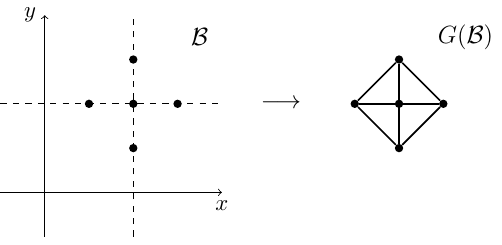}
%
%
%
%
%
\caption{The underlying graph of a metric betweenness structure induced by five points in the Euclidean plane}\label{Fadj}
\end{figure}

The \emph{underlying graph} (or \emph{adjacency graph}) of an almost-metric betweenness structure $\BB = (X,\beta)$ is the graph $G(\BB) = (X, E(\BB))$ where the edges are such pairs of distinct points for which no third point lies between them (see Figure \ref{Fadj}). More formally,
$$E(\BB) =\left\{xz\in\binom{X}{2}:\forall\,y\in X,\,(x\ y\ z)_\BB\Rightarrow\,y = x\vee y = z\right\}.$$

These edges are also called primitive pairs by some authors. The \emph{underlying graph of a metric space} $M$ is $G(M) = G(\BB(M))$. The underlying graph is an important graph-theoretic tool in studying betweenness structures, therefore, it is highly desirable to gain a better understanding of its properties.

The fact that the underlying graph of a metric betweenness structure is connected is part of the mathematical folklore.
Changat et. al. generalized this observation to the path transit function equivalent of almost-metric betweenness structures in \cite{changat2010induced}.

\begin{prp}[Changat et al. \cite{changat2010induced}]\label{Pconn}
The underlying graph of an almost-metric betweenness structure is connected.
\end{prp}


Next, we introduce some definitions for simple and weighted graphs.
A \emph{weighted graph} is a triple $W = (V, E,\omega)$ where $G = (V, E)$ is a simple connected graph and $\omega$ is a positive real-valued function on the set of edges, also called the \emph{edge-weighting} of $W$. We will also denote $V$, $E$, $G$ and $\omega$ by $V(W)$, $E(W)$, $G(W)$ and $\omega_W$, respectively.
The connectivity of $G$ and the positivity of $\omega$ guarantee that a proper weighted graph metric can be defined on $W$ (see below).
We note that any simple connected graph $G = (V, E)$ can be regarded as a weighted graph $W(G)$, with edge-weighting $\omega_G: E\rightarrow\{1\}$, 
thus, any definitions for weighted graphs can be naturally applied to simple connected graphs as well.

Let $W = (V, E,\omega)$ be a weighted graph and $P = v_0 e_1 v_1 e_2\ldots e_\ell v_\ell$ be a path in $W$. The \emph{length} of $P$ is $|P| =\ell$ and the \emph{weight} of $P$ is $\omega(P) =\sum_{i\in [\ell]} \omega(e_i)$. We want to emphasize that the length and the weight of a path are usually different, even though the two coincide in case of simple graphs.
The \emph{metric space induced by the weighted graph $W$} is $M(W) = (V, d_W)$ where $d_W(u, v)$ is the minimum weight of a $u$-$v$ path in $W$, for any $u, v\in V$. Such a path is called a \emph{$u$-$v$ geodesic} in $W$. Note that subpaths of a geodesic are also geodesics. Because of our definition of weighted graphs, $d_W$ is a metric that we call the \emph{weighted graph metric} of $W$. For a simple connected graph $G$, $d_G$ is the usual graph metric defined by the length of shortest paths. We remark that every (finite) metric space $M = (X, d)$ is induced by some weighted graph $W$. For example, take $d$ as the edge-weighting of a complete graph on vertex set $X$.

The \emph{betweenness structure induced by $W$} is the betweenness structure induced by $M(W)$, denoted by $\BB(W)$. For simplicity, we will write $(x\ y\ z)_W$ for $(x\ y\ z)_{\BB(W)}$.
Let us remark that the betweenness relation of a weighted graph $W$ can be described in terms of its geodesics: $(x\ y\ z)_W$ holds if and only if $y$ is on an $x$-$z$ geodesic of $W$.

The weighted graph $Z$ is a \emph{weighted subgraph} of $W$ (in notation, $Z\leq W$) if $G(Z)\leq G(W)$ and $\omega_Z$ is the restriction of $\omega_W$ to $E(Z)$. Let $U$ be a set of vertices that induces a connected subgraph of $G(W)$.
The weighted subgraph of $W$ induced by $U$ is the (uniquely determined) weighted subgraph $W[U]\leq W$ for which $G(W[U]) = G(W)[U]$.
We say that a weighted subgraph $Z\leq W$ is \emph{isometric} if $M(Z)$ is a metric subspace of $M(W)$, i.e. for every $x, y\in V(Z)$, $d_Z(x, y) = d_W(x, y)$. We remark that $\BB(Z)\leq\BB(W)$ in that case.

By graph, we will always mean a finite simple connected graph in the rest of the paper, and we will state explicitly if a graph is weighted. We will use the usual notations $P_n$, $C_n$ and $K_n$ for the (non-weighted) path, cycle and complete graph of order $n$.

A betweenness structure/metric space is
\begin{itemize}
\item \emph{graphic} if it is induced by a graph.
\item \emph{ordered} if it is induced by a path.
\end{itemize}

We will denote the ordered betweenness structure induced by the path $P = x_1x_2\ldots x_n$ by $[x_1, x_2,\ldots, x_n]$. 
We remark that almost-metric betweenness structures are typically not graphic. However, there is an important connection between the underlying graph and graphicity of an almost-metric betweenness structure.

\begin{obs}\label{Ogbg}
For every graph $G$, $G(\BB(G))) = G$. Further, for every al\-most-met\-ric betweenness structure $\BB$, $\BB(G(\BB)) =\BB$ if and only if $\BB$ is graphic.
\end{obs}

\section{Main Results}\label{S2}
In this section, we present the main results of this paper: we introduce and characterize uniquely representable graphs, as well as two related classes of graphs that bound their representations from below and from above (Theorem \ref{Tunique1} and Theorem \ref{Tunique2}). These results can be viewed as a new metric characterization of block graphs and distance-hereditary graphs and --as we point out below-- a generalization of a result of Dress \cite{dress2007category} (Proposition \ref{Pdress}).

\begin{figure}[t]
\centering
\includegraphics{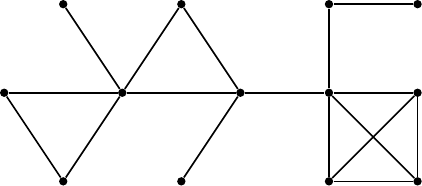}
\caption{A block graph}\label{Fhusimi}
\end{figure}

Graph $G$ is a \emph{block graph} if every block ($2$-connected component) of $G$ is a clique, or equivalently, if every cycle of $G$ induces a complete subgraph (see Figure \ref{Fhusimi}). Block graphs are natural generalizations of trees. They were also called \emph{Husimi-tree}s by some authors, although, that name is not very accurate and refers to another class of graphs today.
$G$ is a \emph{distance-hereditary graph} if all of its connected induced subgraphs are isometric. It can be easily shown that all block graphs are distance-hereditary.

We say that an almost-metric betweenness structure $\BB$ is a \emph{representation} of a graph $G$ if it is the underlying graph of $\BB$. It follows from Observation \ref{Ogbg} that $\BB(G)$ is always a metric representation of $G$.
Below, we investigate the question whether it is possible to infer a betweenness structure from its underlying graph or at least obtain a subset or superset of its betweenness relation. 

\begin{dfn}\label{Dunique}
A graph $G$ is uniquely representable in the metric/\-al\-most-metric sense (ur-m/ur-am) if $\BB(G)$ is the only metric/\-al\-most-metric representation of $G$.
\end{dfn}

We will also consider the following generalizations of Definition \ref{Dunique}. 

\begin{dfn}\label{Dunique2}
A graph $G$
\begin{itemize}
\item \emph{bounds its representations from below} in the metric/almost-metric sense (brb-m/brb-am)  if $\BB(G)\preccurlyeq\BB$ holds for all metric/almost-metric representations $\BB$ of $G$;
\item \emph{bounds its representations from above} in the metric/almost-metric sense (bra-m/bra-am) if $\BB(G)\succcurlyeq\BB$ holds for all metric/almost-metric representations $\BB$ of $G$.
\end{itemize}
\end{dfn}

We denote the corresponding classes of graphs by $\mathcal{C}_{ur}^{m}, \mathcal{C}_{ur}^{am}, \mathcal{C}_{bra}^{m}, \mathcal{C}_{bra}^{am}, \mathcal{C}_{brb}^{m}$ and $\mathcal{C}_{brb}^{am}$, respectively. It is easy to see that the following relations hold between these classes:

\begin{obs}\label{Orel}$ $
\begin{itemize}
\item $\mathcal{C}_{ur}^{m} =\mathcal{C}_{bra}^{m}\cap\mathcal{C}_{brb}^{m}$ and $\mathcal{C}_{ur}^{am} =\mathcal{C}_{bra}^{am}\cap\mathcal{C}_{brb}^{am}$;
\item $\mathcal{C}_{ur}^{am}\subseteq\mathcal{C}_{ur}^{m}$, $\mathcal{C}_{bra}^{am}\subseteq\mathcal{C}_{bra}^{m}$ and $\mathcal{C}_{brb}^{am}\subseteq\mathcal{C}_{brb}^{m}$.
\end{itemize}
\end{obs}

The next two theorems are the main results of this paper. We characterize uniquely representable graphs, as well as graphs that bound their representations from below and from above in both the metric and the almost-metric sense. The proofs can be found in Section \ref{S4}.

\begin{thr}\label{Tunique1}
The following statements are equivalent:
\begin{enumerate}
\item $G\in\mathcal{C}_{bra}^{m}$;
\item $G\in\mathcal{C}_{bra}^{am}$;
\item $G$ is a distance-hereditary graph.
\end{enumerate}
\end{thr}

\begin{thr}\label{Tunique2}
The following statements are equivalent:
\begin{enumerate}
\item $G\in\mathcal{C}_{ur}^{m}$;
\item $G\in\mathcal{C}_{ur}^{am}$;
\item $G\in\mathcal{C}_{brb}^{m}$;
\item $G\in\mathcal{C}_{brb}^{am}$;
\item $G$ is a block graph.
\end{enumerate}
\end{thr}

%


Our motivation for characterizing uniquely representable graphs was two-folded. First, we wanted to have a better understanding on the relationship of betweenness structures and their underlying graphs. We have observed that under certain conditions, a metric betweenness structure can be fully reconstructed from its underlying graph, which is exactly what unique representability means. For example, an interesting remark of Dress implies that trees are uniquely representable in the metric sense.
\begin{prp}[Dress \cite{dress2007category}]\label{Pdress}
Let $\BB$ be a metric betweenness structure such that $T = G(\BB)$ is a tree. Then $\BB$ is induced by $T$.
\end{prp}

We have found Proposition \ref{Pdress} particularly useful in our research on linear betweenness structures.
For example, we were able to obtain a new proof for the main result of \cite{richmond1997metric} by the help of it \cite{fms-lin}. This brings us to our second motivation, which is to find useful generalizations of Proposition \ref{Pdress} that can be applied to a larger set of problems.

Graphs that bound their representation from below and from above arise as the natural generalization of uniquely representable graphs. As Theorem \ref{Tunique1} and Theorem \ref{Tunique2} show, these broader classes of graphs have nice characterizations, too.
We can make two interesting observations on the main results. First, there is no difference between the metric and almost-metric cases: $\mathcal{C}_{ur}^{m} =\mathcal{C}_{ur}^{am}$, $\mathcal{C}_{bra}^{m} =\mathcal{C}_{bra}^{am}$ and $\mathcal{C}_{brb}^{m} =\mathcal{C}_{brb}^{am}$. Second, the results show a significant asymmetry: while $\mathcal{C}_{ur}^{m}$ is a proper subclass of $\mathcal{C}_{bra}^{m}$, it coincides with $\mathcal{C}_{brb}^{m}$ (and the same holds in the almost-metric case).

Finally, we note that in \cite{mulder2008transit}, Prototype problem $3$, Mulder discussed a question involving a kind of unique representability but it is applied to transit functions and convexities instead of betweenness structures and underlying graphs.
Although the concept is quite different from the one discussed here, the example shows that unique representability is a general concept, and investigating related problems for different kinds of combinatorial structures may give rise to interesting new ideas and results.

\section{Preparations}\label{S3}

In this section, we prepare for the proof of Theorem \ref{Tunique1} and Theorem \ref{Tunique2} by
recalling some useful results about weighted graphs, block and distance-hereditary graphs and extending the definition of a geodesic to almost-metric betweenness structures.

\subsection{Block Graphs and Distance-hereditary Graphs}\label{S31}

Block graphs and distance-hereditary graphs are well-understood graph classes and have been characterized in many different ways.
Block graphs were first studied by Harary in \cite{harary1963characterization}. They were characterized in terms of various metric properties (such as the four point condition, the Ptolemaic property and the weakly geodetic property) \cite{kay1965characterization, howorka1979metric}, forbidden induced subgraphs \cite{bandelt1986distance-hereditary}, bulge \cite{howorka1979metric}, k-Steiner intervals \cite{bresar2009steiner}, vertex induced partitions \cite{dress2016characterizing} and betweenness axioms \cite{balakrishnan2015axiomatic}. Further, Le and Tuy gave a characterization and a linear time recognition algorithm for squares of block graphs \cite{le2010square}.

Distance-hereditary graphs were characterized in terms of conditions on paths and cycles by Howorka \cite{howorka1977characterization}, who also showed that Distance-hereditary graphs are perfect.
Bandelt and Mulder gave a good summary on the most important characterizations of distance-hereditary graphs including recursive, metric and forbidden subgraph characterizations \cite{bandelt1986distance-hereditary}, while Morgana and Mulder provided a characterization in terms of axioms imposed on the induced path function \cite{morgana2002induced}.

Below, we list some of the results on which our proof of Theorem \ref{Tunique1} and Theorem \ref{Tunique2} relies on. First, observe the following property of block graphs, which is a straightforward consequence of the definition.

\begin{obs}\label{Ouniquepath}
In a block graph, any two vertices are connected by a unique induced path.
\end{obs}

A \emph{chord} of cycle is an edge that connects two non-consecutive vertices of the cycle. A graph is \emph{chordal} if it does not contain any induced cycles of length at least four, i.e. if it does not contain a chordless cycle. By diamond, we mean a $4$-cycle with exactly one chord.

\begin{prp}[Bandelt, Mulder \cite{bandelt1986distance-hereditary}]\label{Pblock}
Block graphs are exactly the di\-a\-mond\--free chor\-dal graphs.
\end{prp}

\begin{prp}[Howorka \cite{howorka1977characterization}]\label{Pdisthch}
A graph is distance-hereditary if and only if every induced path in it is a geodesic.
\end{prp}

The following statement is an easy consequence of Proposition \ref{Pdisthch}.
\begin{obs}\label{Oindc5}
A distance hereditary graph does not have any induced cycles of length at least $5$.
\end{obs}

\subsection{Tight Weighted Graphs}\label{S32}

We say that an edge $e = xy$ of a weighted graph $W$ is \emph{tight} if $e$ is the unique $x$-$y$ geodesic in $W$, in other words, every $x$-$y$ geodesic in $W$ is of length $1$. A weighted graph is \emph{tight} if all of its edges are tight. Note that all simple graphs are tight. Further, every geodesic in a tight weighted graph $W$ is an induced path in $G(W)$.

\begin{prp}\label{Pweighted}
If $W$ is a tight weighted graph, then $G(W) = G(\BB(W))$.
\end{prp}
\begin{prf}
Let $W = (V, E,\omega)$ be a tight weighted graph. Observe that $G(\BB(W))\leq G(W)$ is always true, since every vertex $y$ on an $x$-$z$ geodesic of $W$ satisfies $(x\ y\ z)_W$. Thus, it suffices to show that $G(W)\leq G(\BB(W))$. Let $e = xz$ be an edge of $G(W)$. If $e\notin E(\BB(W))$, then there exists a vertex $y\in V\backslash\{x, z\}$ such that $(x\ y\ z)_W$ holds, i.e. $d_W(x, z) = d_W(x, y) + d_W(y, z)$. But then, we would get an $x$-$z$ geodesic of length at least $2$ by concatenating an $x$-$y$ and an $y$-$z$ geodesic, in contradiction with the tightness of $e$.
$\square$
\end{prf}

\begin{rmk}
Proposition \ref{Pweighted} can also be reversed: if $G(W) = G(\BB(W))$, then $W$ is a tight weighted graph.
\end{rmk}

\subsection{Geodesics in Betweenness Structures}\label{S33}

In this subsection, we generalize the notion of geodesics to almost-metric betweenness structures and summarize their most important basic properties. Let $\BB$ be an almost-metric betweenness structure on ground set $X$, and let $x, z\in X$ be two points of $\BB$.
An \emph{$x$-$z$ geodesic in $\BB$} is an induced $x$-$z$ path $P$ in $G(\BB)$ such that $\BB\vert_{V(P)}$ is an ordered substructure of $\BB$. The following Proposition lists the most important basic properties of geodesics in betweenness structures.


\begin{prp}\label{Pgeodbetw}
Let $\BB$ be an almost-metric betweenness structure on ground set $X$. Then
\begin{enumerate}
\item for every geodesic $P$ in $\BB$, $\BB\vert_{V(P)} =\BB(P)$;\label{Egeodbetw1} 
\item for every maximal ordered set $Y$ in $\BB$, $G(\BB)[Y]$ is a geodesic in $\BB$;\label{Egeodbetw2} 
\item for every $x, y\in X$, there exists an $x$-$y$ geodesic in $\BB$;\label{Egeodbetw3} 
\item if $\BB$ is induced by a tight weighted graph $W$, then the geodesics of $\BB$ coincide with the geodesics of $W$;\label{Egeodbetw4} 
\item for every $x, y, z\in X$, $(x\ y\ z)_\BB$ holds if and only if there exists an $x$-$z$ geodesic in $\BB$ that contains $y$;\label{Egeodbetw5} 
\item for every almost-metric betweenness structure $\mathcal{C}$ for which $G(\mathcal{C}) = G(\BB)$ holds, $\BB\preccurlyeq\mathcal{C}$ if and only if all geodesics in $\BB$ are geodesics in $\mathcal{C}$.\label{Egeodbetw6} 
\end{enumerate}
\end{prp}
\begin{prf}
Part \ref{Egeodbetw1} is obvious from the definition of geodesics, while Part \ref{Egeodbetw2} follows from the four relations property: if $\BB\vert_Y = [y_1, y_2,\ldots, y_k]$ and $G(\BB)[Y]$ is not a path, then there exists an index $1\leq i< k$ and a point $x\in X\backslash Y$ such that $(y_i\ x\ y_{i+1})_\BB$, which implies that $Y + x$ is ordered by property \ref{Ecoll4}, in contradiction with the maximality of $Y$.


Part \ref{Egeodbetw3} is a simple consequence Part \ref{Egeodbetw2}: take the geodesic induced by a maximal ordered set containing $x$ and $y$ and then take its subpath between $x$ and $y$.

As for Part \ref{Egeodbetw4}, let $W = (V, E,\omega)$ and let $G = G(W)$. Since $W$ is tight, $G(\BB) = G$ by Proposition \ref{Pweighted}, hence, it is enough to show that an induced path $P = y_1y_2\ldots y_\ell$ in $G$ is a geodesic in $W$ if and only if it is a geodesic in $\BB$ (here, we used that all geodesics in $W$ and $\BB$ are induced paths).

\begin{obs}[Polygon-equality]\label{Opolyeq}
Let $\BB$ be a metric betweenness structure induced by a metric space $M = (X, d)$ and $Y =\{y_1, y_2,\ldots, y_\ell\}$ be a non\-empty set of points of $\BB$. Then $d(y_1, y_\ell) =\sum_{i = 1}^{\ell - 1} d(y_i, y_{i + 1})$ if and only if $\BB\vert_Y = [y_1, y_2,\ldots, y_\ell]$.
\end{obs}

Now, the above statement follows from Part \ref{Egeodbetw1} and Observation \ref{Opolyeq}: $P$ is a geodesic in $W$ if and only if $d_W(y_1, y_\ell) = \sum_{i = 1}^{\ell - 1} d_W(y_i, y_{i + 1})$ if and only if $\BB\vert_{V(P)} = [y_1, y_2,\ldots, y_\ell]$ if and only if $P$ is a geodesic in $\BB$.

In order to prove Part \ref{Egeodbetw5}, suppose first that $(x\ y\ z)_\BB$ holds, and let $Y$ be a maximal ordered set in $\BB$ that contains $x, y$ and $z$. By Part \ref{Egeodbetw2}, $P = G(\BB)[Y]$ is a geodesic in $\BB$, and Part \ref{Egeodbetw1} implies $\BB\vert_Y =\BB(P)$ from which $(x\ y\ z)_P$ follows.
Second, suppose that $P$ is an $x$-$z$ geodesic in $\BB$ and $y\in V(P)$, i.e. $(x\ y\ z)_P$ holds. Because of Part \ref{Egeodbetw1}, $\BB\vert_{V(P)} =\BB(P)$, hence, $(x\ y\ z)_\BB$ holds as well.

Finally, Part \ref{Egeodbetw6} is a straightforward consequence of Part \ref{Egeodbetw5}. $\square$
\end{prf}

We conclude this section with a lemma that plays a key role in the proof of the main results.

\begin{lmm}\label{Lnorel}
Let $G$ be a graph that contains an induced path which is not a geodesic in $G$. Then $G\notin\mathcal{C}_{bra}^{m}$ and $G\notin\mathcal{C}_{brb}^{m}$.
\end{lmm}
\begin{prf}
We construct a metric representation $\BB$ of $G$ such that $\BB\not\preccurlyeq\BB(G)$ and $\BB\not\succcurlyeq\BB(G)$.
Let $P$ be an induced $x$-$y$ path in $G$ of length $\ell$ which is not a geodesic. Further, let us define the weighted graph $W$ for which $G(W) = G$ and $\omega_W =\mathbf{1}_{E(G)\backslash E(P)} + \varepsilon\mathbf{1}_{E(P)}$ where $0 <\varepsilon < 1/\ell$ and $\mathbf{1}_A$ denotes the indicator function of the set $A\subseteq E(G)$.

First, we show that $P$ is an $x$-$y$ geodesic in $W$ and it is the only one. Let $Q$ be an $x$-$y$ geodesic in $W$. The weight of $P$ is $\omega(P) =\varepsilon\ell < 1$, hence, $Q$ cannot have any edges from outside of $E(P)$. Therefore, $Q$ must be a subpath of $P$ that connects $x$ and $y$ and so $Q = P$.

Next, we prove that $W$ is tight.
The edges of $P$ are obviously tight because they are the shortest edges of $W$. Now, let $e = uv\in E(G)\backslash E(P)$. If $e$ is not tight, then there exist an $u$-$v$ geodesic $Q'$ of length at least $2$ such that $\omega(Q')\leq\omega(e) = 1$. Again, this can only be true if $Q'$ is a subpath of $P$. But then, $e$ would connect two vertices of $P$, which contradicts the choice of $P$ as an induced path.

Now, let $\BB =\BB(W)$ and let $Q$ be an $x$-$y$ geodesic in $G$. Note that $Q\neq P$ because $P$ was not a geodesic in $G$. On the other hand, we have seen that $P$ is a unique geodesic in $W$, thus, $Q$ cannot be a geodesic in $W$.
Finally, we obtain from Proposition \ref{Pweighted} and Observation \ref{Ogbg} that $G(\BB) = G(W) = G = G(\BB(G))$, so we can conclude by Part \ref{Egeodbetw4} and Part \ref{Egeodbetw6} of Proposition \ref{Pgeodbetw} that $\BB\not\preccurlyeq\BB(G)$ and $\BB\not\succcurlyeq\BB(G)$.
$\square$
\end{prf}


\section{Proof of the Main Results}\label{S4}

\subsection{Proof of Theorem \ref{Tunique1}}\label{S41}
Let $G$ be a fixed graph. We show the following chain of implications: $G\in\mathcal{C}_{bra}^{am}\Rightarrow G\in\mathcal{C}_{bra}^{m}\Rightarrow G$ is a distance-hereditary graph $\Rightarrow G\in\mathcal{C}_{bra}^{am}$.

\begin{step}{1 ($G\in\mathcal{C}_{bra}^{am}\Rightarrow G\in\mathcal{C}_{bra}^{m}$)}
This step follows from Observation \ref{Orel}.
\end{step}

\begin{step}{2 ($G\in\mathcal{C}_{bra}^{m}\Rightarrow G$ \normalfont{is a distance-hereditary graph})}
By Proposition \ref{Pdisthch}, it is enough to show that all induced paths in $G$ are geodesics, which follows from Lemma \ref{Lnorel}.
\end{step}

\begin{step}{3 ($G$ \normalfont{is a distance-hereditary graph} $\Rightarrow G\in\mathcal{C}_{bra}^{am}$)}
Let $\BB$ be an almost-metric representation of $G$. In order to prove $\BB\preccurlyeq\BB(G)$, it suffices to show, by Part \ref{Egeodbetw6} of Proposition \ref{Pgeodbetw}, that every geodesic of $\BB$ is a geodesic of $\BB(G)$. If $P$ is a geodesic of $\BB$, then it is an induced path in $G(\BB) = G$. Therefore, $P$ must be a geodesic in $G$ by Proposition \ref{Pdisthch} and thus $P$ is a geodesic in $\BB(G)$ by Part \ref{Egeodbetw4} of Proposition \ref{Pgeodbetw}. $\square$
\end{step}

%
%

\subsection{Proof of Theorem \ref{Tunique2}}\label{S42}

Let $G = (V, E)$ be a fixed graph. We show the following chain of implications: $G\in\mathcal{C}_{ur}^{m}\Rightarrow G$ is a block graph $\Rightarrow G\in\mathcal{C}_{ur}^{am}\Rightarrow G\in\mathcal{C}_{brb}^{am}\Rightarrow G\in\mathcal{C}_{brb}^{m}\Rightarrow G\in\mathcal{C}_{ur}^{m}$.

\begin{step}{1 ($G\in\mathcal{C}_{ur}^{m}\Rightarrow G$ \normalfont{is a block graph})}
Assume to the contrary that $G$ is not a block graph. It follows from $G\in\mathcal{C}_{ur}^{m}$ that $G$ is distance-hereditary by Theorem \ref{Tunique1}, and we obtain from Proposition \ref{Pblock} and Observation \ref{Oindc5} that $G$ contains an induced subgraph $H$ that is either a diamond or a cycle of length $4$.

\begin{figure}[t]
\centering
\hfill\begin{subfigure}[b]{0.45\textwidth}
\centering
\includegraphics{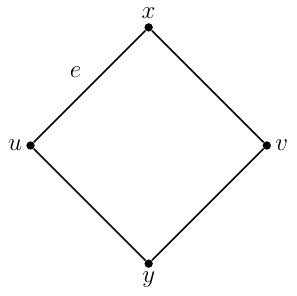}
\caption{$H\simeq C_4$}
\end{subfigure}\hfill
\begin{subfigure}[b]{0.45\textwidth}
\centering
\includegraphics{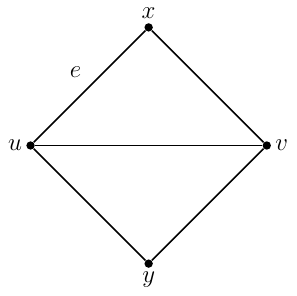}
\caption{$H\simeq D$}
\end{subfigure}\hfill
\caption{Graph $H$ in Step 1 of the proof of Theorem \ref{Tunique2}}\label{Fstep2}
\end{figure}

Let $x$, $y$, $u$ and $v$ denote the vertices of $H$ such that $xy\notin E$, and let $e = xu$, which is an edge of $H$ irrespective of whether $H$ is a diamond or a $4$-cycle (see Figure \ref{Fstep2}). Further, let us define a weighted graph $W = (V, E,\omega)$ where $\omega =\mathbf{1}_E +\frac 1 2\mathbf{1}_{\{e\}}$, and let $\BB =\BB(W)$ and $Z = W[\{x, y, u, v\}]$.

We show that $\BB$ is a representation of $G$ but $\BB\neq\BB(G)$, in contradiction with $G\in\mathcal{C}_{ur}^{m}$. Observe that $W$ is tight because every path of length at least $2$ in $W$ is of weight at least $2$, which is greater than any of the edge weights ($\omega$ is bounded by $\frac 3 2$ from above). Thus, $G(\BB) = G(W) = G$ by Proposition \ref{Pweighted}.
Also notice that $(x\ u\ y)_H$ is true but $(x\ u\ y)_{Z}$ is false, and further, $H$ is an isometric subgraph of $G$ and $Z$ is an isometric weighted subgraph of $W$. This implies that $(x\ u\ y)_G$ holds but $(x\ u\ y)_{W}$ does not, hence, $\BB\neq\BB(G)$.
\end{step}

\begin{step}{2 ($G$ \normalfont{is a block graph} $\Rightarrow G\in\mathcal{C}_{ur}^{am}$)}
Let $G$ be a block graph and $\BB$ be an almost-metric representation of $G$.
By Part \ref{Egeodbetw4} and Part \ref{Egeodbetw6} of Proposition \ref{Pgeodbetw}, it suffices to show that $G$ and $\BB$ have the same geodesics.

Let $x$ and $y$ be any two vertices of $G$, and let $P$ and $P'$ be $x$-$y$ geodesics in $G$ and $\BB$, respectively. Note that such a $P'$ exists by Part \ref{Egeodbetw3} of Proposition \ref{Pgeodbetw}.
Further, since both $P$ and $P'$ are induced $x$-$y$ paths in $G$, Observation \ref{Ouniquepath} yields $P = P'$, which completes the proof.
\end{step}

\begin{step}{3 ($G\in\mathcal{C}_{ur}^{am}\Rightarrow G\in\mathcal{C}_{brb}^{am}$)}
This step follows from Observation \ref{Orel}.
\end{step}

\begin{step}{4 ($G\in\mathcal{C}_{brb}^{am}\Rightarrow G\in\mathcal{C}_{brb}^{m}$)}
Again, this step follows from Observation \ref{Orel}.
\end{step}

\begin{step}{5 ($G\in\mathcal{C}_{brb}^{m}\Rightarrow G\in\mathcal{C}_{ur}^{m}$)}
If $G\in\mathcal{C}_{brb}^{m}$, then because of Lemma \ref{Lnorel}, every induced path in $G$ must be a geodesic, hence, $G$ is distance hereditary by Proposition \ref{Pdisthch}. Finally, we obtain from Theorem \ref{Tunique1} that $G\in\mathcal{C}_{bra}^{m}$, therefore, $G\in\mathcal{C}_{ur}^{m}$ by Observation \ref{Orel}. $\square$
\end{step}

\section{Conclusion}\label{S5}
Motivated by our observations on finite metric spaces, we have defined the class of uniquely representable graphs and the class of graphs that bound their representations from below and from above. 
These graph-classes have been characterized in Theorem \ref{Tunique1} and Theorem \ref{Tunique2} in both the metric and the almost-metric sense. In particular, we have shown that the uniquely representable graphs are exactly the block graphs and pointed out that this result generalizes an interesting remark of Dress about trees (Proposition \ref{Pdress}).

Lastly, we would like to mention two open problems on graph representability that can be potential subjects of future research.
By definition, uniquely representable graphs have the minimum number of representations. One may also be interested in the other extreme.

\begin{prb}
What is the maximum number of metric/almost-metric representations of a graph on $n$ vertices and which graphs realize that maximum?
\end{prb}
We conjecture the following.

\begin{cnj}
The number of metric representations of a graph of order $n$ is maximized by the balanced complete bipartite graph of order $n$.
\end{cnj}

Note that the balanced complete bipartite graph $K_{\left\lfloor\frac{n}{2}\right\rfloor,\left\lceil\frac{n}{2}\right\rceil}$ on $n$ vertices have at least $2^{\left\lfloor\frac{n}{2}\right\rfloor\left\lceil\frac{n}{2}\right\rceil - n + 1}$ metric representations. Namely, pick one vertex from both classes of $K_{\left\lfloor\frac{n}{2}\right\rfloor,\left\lceil\frac{n}{2}\right\rceil}$, and set the weight of the edges adjacent to these vertices to $1$. For all other edges, choose a weight arbitrarily from the set $\{1, 2\}$. It can be easily seen that different weighted graphs obtained in this way induce different metric representations of $K_{\left\lfloor\frac{n}{2}\right\rfloor,\left\lceil\frac{n}{2}\right\rceil}$.

We have weakened the property of unique representability by only requiring $\BB(G)$ to be a smallest/largest element in the poset of representations of $G$, where elements are ordered according to relation $\preccurlyeq$. We can further weaken this condition by saying that $\BB(G)$ is a minimal/maximal element in the poset of representations.

\begin{prb}
Characterize graphs $G$ for which $\BB(G)$ is a minimal/maximal element in the poset of its metric/almost-metric representations.
\end{prb}

\section*{Acknowledgment}
The author is grateful to Pierre Aboulker for sharing his thoughts on the results presented here, as well as on other, related results.

This work was supported by the National Research, Development and Innovation Office -- NKFIH, No. 108947.

\bibliographystyle{elsarticle-num}
\bibliography{fms_bib}

\end{document}